\documentclass[12pt]{article} 
\usepackage{graphicx}
\newcommand{\be}{\begin{equation}}
\newcommand{\ee}{\end{equation}}
\newcommand{\ba}{\begin{array}}
\newcommand{\ea}{\end{array}}

\newcommand{\RR}{\mathbb R}

\newcommand{\bea}{\begin{eqnarray*}}
\newcommand{\eea}{\end{eqnarray*}}
\newcommand{\bean}{\begin{eqnarray}}
\newcommand{\eean}{\end{eqnarray}}

\def\ds{\displaystyle}
\def\nm{\noalign{\medskip}}

\newtheorem{lemma}{Lemma}[section]

\newtheorem{theorem}{Theorem}[section]

\usepackage{amsmath, amssymb}
\usepackage{latexsym}

\def\Box{\leavevmode\vbox{\hrule
     \hbox{\vrule\kern5pt\vbox{\kern5pt}%
           \vrule}\hrule}}

\numberwithin{equation}{section}
\newcommand{\Scal}{\mathcal{S}}
\newcommand{\Dcal}{\mathcal{D}}
\newcommand{\Kcal}{\mathcal{K}}

\newcommand{\Rcal}{\mathcal{R}}
\newcommand{\Lcal}{\mathcal {L}}

\begin{document}
\title{Dynamic Linear Equations: Small Perturbations of an
Interface }


\author{Christian Daveau \thanks{
 Department of Mathematics, CNRS (UMR 8088),
 University of Cergy-Pontoise,
 2 avenue Adolphe Chauvin,
 95302 Cergy-Pontoise Cedex, France. (Email: christian.daveau@u-cergy.fr)}
\and Abdessatar Khelifi \thanks{
D\'epartement de Math\'ematiques, Universit\'e des Sciences de
Carthage, Bizerte, Tunisie. (Email: abdessatar.khelifi@fsb.rnu.tn)}
}

\maketitle \abstract{We consider solutions to the initial boundary value problems for both heat and wave equations in two dimensions. In this paper, we aim to advance the development of high-order terms in the asymptotic expansions of current boundary perturbations caused by small perturbations of the shape of a $C^2$-boundary inhomogeneity. Our derivations are rigorously derived by the field expansion (FE) method (formal derivation) and proved by layer potential techniques. It extends those already derived for small-volume inhomogeneities to develop efficient algorithms for determining certain properties of the shape of an inhomogeneity based on boundary measurements.}\\

.


\section{Introduction}
Let $\Omega$ be a bounded domain in $\mathbb{R}^2$, with a connected $C^2$-boundary $\partial\Omega$. We suppose that $\Omega$ contains an inhomogeneity $B$ well separated from the boundary $\partial\Omega$. Let $1$ denote the conductivity of the background medium, while $k$ denotes the constant conductivity of the inhomogeneity $B$.
In this project, we suppose that $B_\alpha$ is an $\alpha$-perturbation of $B$, i.e., there is a function $h \in C^2(\partial B)$ such that
$$
\partial B_\alpha:=\{\tilde{x}=x+\alpha h(x)\nu(x):=\Psi_\alpha(x),~ x \in \partial B \}, \mbox{ where }0 < \alpha << 1.$$
In the presence of the above $\alpha$-perturbation, we consider the initial boundary value problem for the perturbed heat equation:
\begin{equation}\label{u-heat-pert}
\left\{
  \begin{array}{ll}
    \partial_t U_\alpha-\nabla\cdot \Big(\chi_{\Omega\setminus \overline{B}_\alpha}+k\chi_{B_\alpha}
\Big)\nabla U_\alpha = 0 & \mbox{ in } \Omega\times (0,T), \\
    U_\alpha(x,0)=\phi(x)&\mbox{ in } \Omega,\\
   \frac{\partial U_\alpha}{\partial\nu}=f &\mbox{ on }\partial\Omega\times(0,T),\\
    U_\alpha=f&\mbox{ on }\partial\Omega\times(0,T).
  \end{array}
\right.
\end{equation}

and so, the perturbed scalar wave problem is given by:
\begin{equation}\label{u-wave-pert}
\left\{
  \begin{array}{ll}
    \partial_t^2 u_\alpha-\nabla\cdot \Big(\chi_{\Omega\setminus \overline{B}_\alpha}+k\chi_{B_\alpha}
\Big)\nabla u_\alpha = 0&\mbox{ in } \Omega\times (0,T), \\
    u_\alpha(x,0)=\varphi(x), \quad \partial_t u_\alpha(x,0)=\psi(x)&\mbox{ in } \Omega,\\
    \frac{\partial u_\alpha}{\partial\nu}=f&\mbox{ on }\partial\Omega\times(0,T), \\
    u_\alpha=f&\mbox{ on }\partial\Omega\times(0,T).
  \end{array}
\right.
\end{equation}
Here $\chi_B$ is the indicator function of $B$, $T > 0$ is a final observation time, $\nu$ denotes the unit outward normal to the domain $\Omega$
, $\phi(x)$ the initial data, $g(x, t)$ is a suitably prescribed function, $f(x, t)$ is the flux, and $\varphi(x)$ is the initial data.\\

In the absence of perturbation, we consider the initial boundary value problem for the scalar heat equation

\begin{equation}\label{main-heat0}
\left\{
  \begin{array}{ll}
    \partial_t U_0-\nabla\cdot \Big(\chi_{\Omega\setminus \overline{B}}+k\chi_{B}
\Big)\nabla U_0 = 0 & \mbox{ in } \Omega\times (0,T), \\
    U_0(x,0)=\phi(x)&\mbox{ in } \Omega,\\
   \frac{\partial U_0}{\partial\nu}=f &\mbox{ on }\partial\Omega\times(0,T),\\
    U_0=f&\mbox{ on }\partial\Omega\times(0,T),
  \end{array}
\right.
\end{equation}
while the unperturbed scalar wave problem is given by
\begin{equation}\label{main-wave0}
\left\{
  \begin{array}{ll}
    \partial_t^2 u_0-\nabla\cdot \Big(\chi_{\Omega\setminus \overline{B}}+k\chi_{B}
\Big)\nabla u_0 = 0&\mbox{ in } \Omega\times (0,T), \\
    u_0(x,0)=\varphi(x), \quad \partial_t u_0(x,0)=\psi(x)&\mbox{ in } \Omega,\\
    \frac{\partial u_0}{\partial\nu}=f&\mbox{ on }\partial\Omega\times(0,T), \\
    u_0=f&\mbox{ on }\partial\Omega\times(0,T).
  \end{array}
\right.
\end{equation}

Our goal in this project is to write explicit high order asymptotic formulae for both heat and wave perturbed problems and to reconstruct the
shape deformation $ h$ from dynamical measurements of boundary perturbations $U_\alpha-U$ for heat conduction problem (resp.
$u_\alpha-u$ for wave equation) on $\partial\Omega$.\\

Our dynamic identification procedure can be regarded as a first attempt to generalize that
of [1,2,3,4,5,6] to the time-dependent case. It yields simple and effective computational identification algorithms.\\
We will perform numerical experiments to test the viability of the proposed algorithms. We believe that our results are ambitious tools for
determining the locations and/or shapes of small inhomogeneities by taking dynamical measurements.\\


In a similar approach, the first order term in voltage potential boundary perturbations due to thin interfaces with non-extreme conductivity has been formally derived and justified in \cite{BMV} and \cite{BFV}. Investigations have been conducted on the behavior of voltage potentials in the presence of thin interfaces of high or low conductivity. Various aspects related to the effective conductivity of composite materials with thin interfaces of either low or high conductivity have been studied. To the best of our knowledge, this
is the first work to rigorously derive leading order terms in the
asymptotic expansions  of $u_\alpha  |_{\partial \Omega}$ when  $\alpha  \to 0$ and
no uniform thickness, that is, $h$ is not constant. Proofs of our
asymptotic expansions in the three cases which are mentioned above
are radically different from those in \cite{BFV,
palencia1,palencia2} and  based on a general unified layer potential
technique. Our procedure can be continued and all the terms in the
expansion of $u_\alpha  |_{\partial \Omega}$ can be obtained in the same fashion.
Compared to several previous works, our method seems simpler to implement,
 and it has the big advantage to allow a generalization to $n$ dimensions
  conductivity inclusions which requires precise calculations of the multi-dimensional approaches.
   In the three-dimensional case, our derivation generalizes the results already presented at the beginning of the paper.
 Thus, we assume that $\partial D_{\alpha }$ has a regular parametrization as a revolution or a general surface.\\

The asymptotic expansions in this paper can be used to design
algorithms to identify certain properties of their thin like
location and thickness.

Our general approach can be extended to other equations such as the
anisotropic conductivity problem, Stokes, the Maxwell and the
Lam{\'e} systems.

This paper is organized as follows. In the next section we state the
main results and investigate the particular case where the leading order
terms are expressed explicitly in this case. The remaining part of
the paper is dedicated to proofs of the main results stated in
section $4$. Then, in  section $3$,  we review some basic facts
about the layer potentials and give a representation formula for the
heat potential associated with the thin $D_\alpha
\backslash \overline{D}$. In section $4$, we derive asymptotic
expansions of layers potentials. In section $5$, we derive
asymptotic expansions of heat and wave functions caused by
the presence $D_\alpha  \backslash \overline{D}$, considering separately
these cases. Finally in section $6$, our main results are
leaded algorithm reconstruction methods.


\section{Formal derivations: field expansion (FE) method}
Let $a, b \in \RR,$ with $a<b$, and let $X(t): [a,b]\to \RR^2$ be
the arclength parametrization of $\partial D$, namely, $X$ is an
$\mathcal{C}^2$-function satisfying $|X'(t)|= 1$ for all $t \in
[a,b]$ and
 $$
 \partial D:=\{x=X(t),  t\in [a,b]\}.
 $$
Then the outward unit normal to $D$, $\nu(x)$, is given by
$\nu(x)=R_{-\frac{\pi}{2}}X'(t)$, where $R_{-\frac{\pi}{2}}$ is the
rotation by $-\frac{\pi}/{2}$,  the tangential vector at $x$, $T(x) =
X'(t)$, and $X'(t)\perp X''(t)$. Set the curvature $\tau(x)$ to be
defined by
 $$
 X''(t)=\tau (x)\nu(x).
 $$
We will sometimes use $h(t)$ for $h(X(t))$ and $h'(t)$ for the
tangential derivative of $h(x)$.

Then, $\tilde{X}(t)=X(t)+\alpha  h(t)\nu(x)=X(t)+\alpha
 h(t)R_{-\frac{\pi}{2}}X'(t)$ is a parametrization of ${D}_\alpha$. By
$\tilde{\nu} (\tilde{x})$, we denote the outward unit normal to
$D_\alpha$ at $\tilde{x}$. Then, it is proved in \cite{AKLZ1} that

 \begin{equation} \label{eq:nuexpan}
 \tilde{\nu}(\tilde{x}) = \nu(x) - \alpha h'(t) T(x) + O(\alpha^2), \quad \mbox{uniformly } x \in \partial
 D.
 \end{equation}
The following uniformly convergent expansion for the length element
$d{\sigma_\alpha }(\tilde{y})$ can be obtained easily (see
\cite{AKLZ1}).
 \begin{equation} \label{sigexp}
 d{\sigma_\alpha }(\tilde{y}) = (1 -\alpha \tau(y)h(y)+ O(\alpha ^2))
 d\sigma(y) .
 \end{equation}

Let $\Psi_\alpha $ be the diffeomorphism from $\partial D$ onto
$\partial D_\alpha $ given by $\Psi_\alpha (x) =  x +\alpha  h(t)\nu(x),$
where $x = X(t)$.

\section{Representation formula} \label{Sec:represenation_formula}
Here we consider initial boundary value problems for the heat equation using the thermal potential representation for the solution. Depending on the choice of representation, we are led to a solution of the various boundary integral equations.\\

Here we review some basic facts about the layer potentials which
will allow us to express the solution $u_\alpha $ in Lemma
\ref{lem:rep_formula} by means of layer potentials which are
solution to a particular system of integral equations.\\
We introduce the fundamental solution of the heat equation
$$
\Gamma(x,t)=\frac{H(t)}{4\pi t}\exp(-|x|^2/4t),
$$
where $H ( t )$ is the Heaviside function. The classical single-layer and
double-layer heat potentials of the density functions
$\phi$, $\psi$ on $D$ are defined by,
\begin{align}
\Scal_D \phi (x,t) &:= \int_{0}^{t}\int_{\partial  D} \Gamma (x-y,t-\tau) \phi (y,\tau) d \sigma (y)d\tau
, \quad (x,t) \in Q_D^T\cup Q_{T,D}^c, \label{single}\\ \nm  \Dcal_D \psi (x,t) &:=
 \int_{0}^{t}\int_{\partial  D} \frac{\partial }{\partial \nu_y}\Gamma (x-y,t-\tau) \psi (y,\tau) d \sigma (y)d\tau , \quad x
\in (x,t) \in Q_D^T\cup Q_{T,D}^c, \label{double}
\end{align}
where $Q_D^T=D\times(0,T)$, $Q_{T,D}^c=D^c\times(0,T)$ and $D^c:=\RR^2\backslash\overline{D}$. For $\phi$ and $\psi$ sufficiently smooth we consider the heat potential
\begin{equation}\label{heat-repre}
  U=\Scal_D \phi- \Dcal_D \psi \quad \mbox{in } Q_D^T.
\end{equation}
For a function $u$ defined on $\RR^2 \setminus \partial  D$, we denote
 $$
 \frac{\partial u}{\partial \nu}\Big|_\pm (x) := \lim_{t \to 0^+} \langle \nabla u(x \pm t
\nu_x), \nu_x \rangle, \quad x \in \partial  D, $$ if the limit exists. The
notation $u|_{\pm}$ is understood likewise.

The proof of the following trace formula can be found in
\cite{Folland76}:
\begin{align}
\frac{\partial }{\partial \nu^\pm} \Scal_D \phi(x) &= \biggl(\pm \frac{1}{2} I +
\Kcal_D^*\biggr) \phi (x), \quad x \in \partial  D \label{singlejump}\\
(\Dcal_D \phi)|_{\pm} & = \biggl(\mp \frac{1}{2} I + \Kcal_D\biggr)
\phi (x), \quad x \in \partial  D, \label{doublejump}
\end{align}
where $$ \Kcal _D \phi (x) = \frac{1}{2\pi}  \int_{\partial  D} \frac{\langle
y-x, \nu(y) \rangle }{|x-y|^2} \phi(y) d \sigma(y) $$ and $\Kcal_D^*$ is
the $L^2$-adjoint of $\Kcal_D$. The identity \eqref{doublejump}
implies that $\Dcal_D \phi$ has a jump across $\partial  D$. But
$\frac{\partial (\Dcal_D \phi)}{\partial \nu}$ does not have a jump across $\partial  D$. In
fact, if $\phi \in \mathcal{C}^{1,\eta}(\partial  D)$ for some $\eta
>0$, then for $x\in \partial  D$
 \begin{align*}
 \frac{\partial (\Dcal_D \phi)}{\partial \nu} \Big |_{\pm}(x) &=
 - \frac{1}{2\pi} \int_{\partial  D} \frac{\langle \nu(x), \nu(y)\rangle }{|x-y|^2}
  \big [\phi(y)-\phi(x) \big]
 d\sigma(y)\\
 & \quad + \frac{1}{\pi} \int_{\partial  D}\frac{\langle   x-y,\nu(x)\rangle  \langle   x-y,\nu(y)\rangle  }{|x-y|^4}
  \big[\phi(y)-\phi(x) \big] d\sigma(y),
 \end{align*}
and hence
 \begin{equation} \label{jumpd}
 \frac{\partial (\Dcal_D \phi)}{\partial \nu} \Big |_{+}(x)= \frac{\partial (\Dcal_D \phi)}{\partial \nu} \Big
 |_{-}(x), \quad x \in \partial  D.
 \end{equation}

Let $L^2_0 (\partial  D) : = \{ f \in L^2(\partial  D) : \int_{\partial  D} f d\sigma=0
\}$. The following results are of importance to us. For proofs see
\cite{Folland76}.
\begin{lemma}\label{lem:Verchota}
The operator $\lambda I - \Kcal_D^*$ is invertible on $L^2_0 (\partial  D)$
if\/ $|\lambda| \geq \frac{1}{2}$, and for $\lambda \in (-\infty,
-\frac{1}{2}] \cup (\frac{1}{2}, + \infty )$, $\lambda I -
\Kcal_D^*$ is invertible on $L^2 (\partial  D)$.
\end{lemma}

Let $\lambda ={(k+1)}/({2(k-1)})$. Denote by $\Scal_{D_\alpha  },
\Dcal_{D_\alpha  }, \Kcal_{D_\alpha  }$, and $\Kcal_{D_\alpha  }^*$ the layer
potentials on $D_\alpha  $. The following lemma can be proved in the
exactly same manner as the representation formula in \cite{KS96}. Based on the analysis above, we are now in the position to present the integral representation for the perturbation filed.

\begin{lemma}\label{lem:rep_formula}
The solution $u_\alpha  $ of the problem \eqref{u-wave-pert} can be
represented  as follows:
\begin{equation}\label{300}
u_\alpha  (x)=H_\alpha  (x)+\Scal_{ D_\alpha  }\tilde{\phi_\alpha  } (x)+\Scal_{
D}\psi_\alpha   (x),
\end{equation}
where the harmonic part $H_\alpha  $ of $u_\alpha  $ is given by $H_\alpha  (x)=-
\Scal_{ \Omega  }(g)(x)+ \Dcal_{\Omega }(u_\alpha  |_{\partial \Omega })(x)$, $x\in \Omega $,
and the pair $(\tilde{\phi_\alpha  }, \psi_\alpha  )\in L_0^2(\partial
{D_\alpha  })\times L_0^2(\partial  D)$ is the unique solution of the following
system of integral equations:
\begin{equation} \label{eq:syst1}  \left\{
\begin{array}{l} \ds (\lambda I - \Kcal^*_{ D_\alpha  }
)\tilde{\phi_\alpha  }(\tilde{x})-\frac{\partial  \Scal_D\psi_\alpha  }{\partial
\tilde{\nu}}(\tilde{x}) =\frac{\partial  H_\alpha  }{\partial  \tilde{\nu}}(\tilde{x}),
\quad \tilde{x} \in \partial   D_\alpha  ,\\ \nm \ds (\lambda I + \Kcal^*_{D}
)\psi_\alpha  (x)+\frac{\partial  \Scal_{D_\alpha  }\tilde{\phi_\alpha  }}{\partial  \nu}(x)
=-\frac{\partial  H_\alpha  }{\partial  \nu}(x), \quad  x\in  \partial  D.
\end{array} \right.\end{equation}
\end{lemma}
By integration by parts, one can check that
$$
\tilde{\phi_\alpha  }=(k-1)\frac{\partial  u_\alpha  }{\partial  \nu}\big|_{-}\quad \mbox{on
} \partial  D_\alpha  ,\quad \quad \tilde{\psi_\alpha  }=(1-k)\frac{\partial  u_\alpha  }{\partial
\nu}\big|_{+}\quad \mbox{on } \partial  D.
$$
In this next section, we will provide  asymptotic expansions of the layer potentials,  which are appeared  in the integral equation
system \eqref{eq:syst1} with $(\tilde{\phi_\alpha  }, \psi_\alpha  )\in
C^{1,\eta}(\partial  D_\alpha  )\times C^{1,\eta}(\partial  D)$. These asymptotic
expansions will help us to derive leading terms in the asymptotic
expansions of the solution $u_\alpha $. To do this, we may use the
following estimates
\begin{equation}\label{estmation}
\displaystyle \left\|\frac{\partial
u_{\alpha }}{\partial\nu}|_{\pm}\right\|_{C^{1,\eta}(\partial
D)},\left\|\frac{\partial
u_{\alpha }}{\partial\nu}|_{\pm}\right\|_{C^{1,\eta}(\partial
D_{\alpha })}\leq C,
\end{equation}
for a some positive constant $C$. We think those estimates can be
derived by following the ideas developed in \cite{LV} with
appropriate changes.

\section{Asymptotic expansions of layer potentials}
\label{Sec:Asymptotic_expansions}
Let $a, b \in \RR,$ with $a<b$, and let $X(t): [a,b]\to \RR^2$ be
the arclength parametrization of $\partial D$, namely, $X$ is an
$\mathcal{C}^2$-function satisfying $|X'(t)|= 1$ for all $t \in
[a,b]$ and
 $$
 \partial D:=\{x=X(t),  t\in [a,b]\}.
 $$
Then the outward unit normal to $D$, $\nu(x)$, is given by
$\nu(x)=R_{-\frac{\pi}{2}}X'(t)$, where $R_{-\frac{\pi}{2}}$ is the
rotation by $-{\pi}/{2}$,  the tangential vector at $x$, $T(x) =
X'(t)$, and $X'(t)\perp X''(t)$. Set the curvature $\tau(x)$ to be
defined by
 $$
 X''(t)=\tau (x)\nu(x).
 $$
We will sometimes use $h(t)$ for $h(X(t))$ and $h'(t)$ for the
tangential derivative of $h(x)$.

Then, $\tilde{X}(t)=X(t)+\alpha  h(t)\nu(x)=X(t)+\alpha
h(t)R_{-\frac{\pi}{2}}X'(t)$ is a parametrization of ${D}_\alpha$. By
$\tilde{\nu} (\tilde{x})$, we denote the outward unit normal to
$D_\alpha $ at $\tilde{x}$. Then, it is proved in \cite{AKLZ1} that

 \begin{equation} \label{eq:nuexpan}
 \tilde{\nu}(\tilde{x}) = \nu(x) - \alpha  h'(t) T(x) + O(\alpha ^2), \quad \mbox{uniformly } x \in \partial
 D.
 \end{equation}
The following uniformly convergent expansion for the length element
$d{\sigma_\alpha }(\tilde{y})$ can be obtained easily (see
\cite{AKLZ1}.
 \begin{equation} \label{sigexp}
 d{\sigma_\alpha }(\tilde{y}) = (1 -\alpha  \tau(y)h(y)+ O(\alpha ^2))
 d\sigma(y) .
 \end{equation}

Let $\Psi_\alpha $ be the diffeomorphism from $\partial D$ onto
$\partial D_\alpha $ given by $\Psi_\alpha (x) =  x +\alpha  h(t)\nu(x),$
where $x = X(t)$. The asymptotic expansion of $\Kcal^*_{D_\alpha }$
is given by the following lemma, see \cite{AKLZ1} for a detailed
proof.
\begin{lemma} \label{lemkd}
Let $N\in \mathbb{N}$. There exists $C$ depending only on $N$, the
$\mathcal{C}^2$-norm of $\partial D$ and the $\mathcal{C}^1$-norm of
$h$ such that for any $\tilde{\phi} \in L^2(\partial D_\alpha )$,
 $$
 || (\Kcal^*_{D_\alpha }
 \tilde{\phi}) \circ \Psi_\alpha  - \Kcal_D^*  \phi -
 \sum_{n=0}^{N-1} \alpha ^{n+1} \Kcal_D^{(n+1)}  \phi
 ||_{L^2(\partial D)} \leq C \alpha ^{N+1} || \phi
 ||_{L^2(\partial D)},
 $$
where $\phi := \tilde{\phi} \circ \Psi_\alpha $. The operators
$\Kcal_D^{(n)}$, $n>0$, are compact on $L^2(\partial \Omega)$ and
can be explicitly computed.
\end{lemma}

In particular, $\Kcal_D^{(1)}$ is given by
 \begin{align*}
 \Kcal_{D}^{(1)}\phi(x) &=\tau(x) h(x)\Kcal_{D}^*(
 \phi)(x)-\Kcal_{D}^*(\tau h
 \phi)(x) \\
 & \quad -\frac{1}{\pi}\int_{\partial
 D}\frac{\langle  x-y,\nu(x)\rangle   \langle  x-y,h(t)\nu(x)-h(s)\nu(y)\rangle   }{|x-y|^4}\phi(y)d\sigma(y)\\
 & \quad +\frac{1}{2\pi}\int_{\partial  D}\frac{\langle  h(t)\nu(x)-h(s)\nu(y),\nu(x)\rangle
 }{|x-y|^2}\phi(y)d\sigma(y) \\
 & \quad -\frac{1}{2\pi}\int_{\partial D}\frac{\langle
 x-y,h'(t)T(x)+ \tau(x)h(x)\nu(x)\rangle }{|x-y|^2}\phi(y)d\sigma(y).
 \end{align*}
 Since
 $$
|x-y|^2= |\langle x-y,\nu (x)\rangle  |^2+ |\langle  x-y,T(x)\rangle|^2.
 $$
If $\phi \in \mathcal{C}^{1, \eta}(\partial   D)$, then we get
 \begin{align}
 \Kcal_{D}^{(1)}\phi(x) & =\tau(x) h(x)\Kcal_{D}^*(
 \phi)(x)-\Kcal_{D}^*(\tau h \phi)(x) \nonumber \\
&\quad - \frac{h(x)}{2\pi} \Big[ \int_{\partial D}\frac{1}{|x-y|^2}
\phi(y)d\sigma(y) -
 2 \int_{\partial D} \frac{\langle  x-y,T(x)\rangle  ^2}{|x-y|^4}\phi(y)d\sigma(y) \Big] \nonumber \\
& \quad-\frac{1}{2\pi}\int_{\partial D}\frac{\langle
 x-y,h'(t)T(x)+\tau (x)h(x) \nu(x)\rangle   }{|x-y|^2}\phi(y)d\sigma(y) \nonumber \\
 & \quad - \frac{1}{2\pi} \Big[ \int_{\partial D} \frac{\langle  \nu(x), \nu(y)\rangle   }{|x-y|^2}
 h(y)\phi(y)d\sigma(y) \nonumber \\
 & \quad - 2 \int_{\partial D}\frac{\langle  x-y,\nu(x)\rangle   \langle  x-y,\nu(y)\rangle   }{|x-y|^4}h(y) \phi(y)d\sigma(y) \Big] \nonumber \\
  & =\tau(x) h(x)\Kcal_{D}^*(
 \phi)(x)-\Kcal_{D}^*(\tau h \phi)(x)\nonumber \\
 & \quad -\frac{d}{dt}\big(h \frac{d}{dt}\Scal_D(\phi)\big)(x)
 + \frac{\partial (\Dcal_D(h\phi))}{\partial \nu}(x),
\label{kone}
 \end{align}
where $\frac{d}{dt}$ is the tangential derivative in the direction
of $T(x)=X'(t)$.

Now we are going to derive an asymptotic expansion of
$\frac{\partial (\Scal_D\psi)}{\partial \tilde{\nu}}(\tilde{x})$ for $\psi \in
\mathcal{C}^{1, \eta}(\partial  D)$, $\tilde{x} \in \partial D_\alpha $. Since $\partial  D$
is $\mathcal{C}^{2,\eta}$, $\Scal_D\psi$ is
$\mathcal{C}^{2,\eta}(\RR^2 \setminus D)$. Thus, it follows from
\eqref{eq:nuexpan} and Taylor expansion that, for $x\in \partial D$,
 \begin{align}\nonumber
 \frac{\partial (\Scal_D\psi)}{\partial \tilde{\nu}}(\tilde{x}) &= \tilde{\nu}(\tilde{x}) \cdot
 \nabla \Scal_D\psi(\tilde{x}) \\\nonumber
 &= \Big( \nu(x) - \alpha  h'(t) T(x) + O(\alpha ^2) \Big) \\\nonumber
 & \quad \cdot \Big( \nabla
 \Scal_D\psi \big|_{+}(x) + \alpha  h(t) \sum_{j=1}^2 \partial_j(\nabla
 \Scal_D\psi(x)) \nu_j(x)\big|_{+} + O(\alpha ^{1+\eta}) \Big) \\\nonumber
 &= \frac{\partial (\Scal_D\psi)}{\partial \nu} \bigg|_+ (x) + \alpha  \Big[ -h'(t)
 \frac{\partial (\Scal_D\psi)}{\partial T}(x)+ h(t) \langle  D^2 (\Scal_D\psi)(x) \nu(x),
 \nu(x) \rangle   \Big |_{+}\Big] \\\label{eq:asym_exp_S_1}
 & \quad + O(\alpha ^{1+\eta}),
 \end{align}
where $D^2(\Scal_D\psi)$ denotes the Hessian of $\Scal_D\psi$ and
$O(\alpha^{1+\eta})$ term is bounded by $C \alpha^{1+\eta} \| \psi
\|_{\mathcal{C}^{1,\eta}(\partial D)}$.\\

Recalling the jump relation \eqref{singlejump} and using the
Laplacian expression in local coordinates
 $$
\Delta=\frac{\partial^2}{\partial \nu^2} -\tau \frac{\partial }{\partial
\nu}+\Big(\frac{d}{dt}\Big)^2,
 $$
we obtain the following identity
 $$
 \langle  D^2 (\Scal_D\psi)(x) \nu(x),
 \nu(x) \rangle   |_{\pm}= \tau (\pm \frac{1}{2}I+\Kcal_{D}^*)(\psi)(x)-
  \Big(\frac{d}{dt}\Big)^2(\Scal_D\psi)(x),\quad x\in \partial D,
 $$
which, after substitution into \eqref{eq:asym_exp_S_1} leads to the
following expansion of the single layer potential
\begin{equation} \label{eq:stilde}
\frac{\partial (\Scal_D\psi)}{\partial \tilde{\nu}}(\tilde{x}) = \Big(\frac{1}{2} I
 + \Kcal_D^* \Big) \psi(x) + \alpha  \Rcal_D \psi (x) + O(\alpha^{1+\eta}),
\end{equation}
where
\begin{align} \label{rform}
 \Rcal_D \psi (x) =  \tau (x) h(x)(\frac{1}{2}I+\Kcal_{D}^*)(\psi)(x) -\frac{d}{dt}\big(h
 \frac{d}{dt}\Scal_D(\psi)\big)(x).
\end{align}
 It is easy
to see that $$
 \| \Rcal_D \psi \|_{C^1(\partial D)} \leq C \| \psi \|_{_{C^1(\partial D)}}.$$
We now expand $\frac{\partial (\Scal_{D_\alpha }) \tilde{\phi}}{\partial \nu}(x)$
for $x \in \partial D$ when $\tilde{\phi} \in \mathcal{C}^{1, \eta}(\partial
D_\alpha )$. Let $f$ be a $\mathcal{C}^{1, \eta}$ function on $\partial  D$ and
let $u$ be the solution to $\Delta u=0$ in $D$ satisfying $u=f$ on
$\partial D$. Then, we get
 \begin{align*}
 \int_{\partial  D} \frac{\partial  (\Scal_{D_\alpha } \tilde{\phi})}{\partial  \nu}(x) f(x)
 d\sigma(x)
 & = \int_{\partial D} \Scal_{D_\alpha } \tilde{\phi}(x) \frac{\partial u}{\partial \nu}(x)
 d\sigma(x) \\
 & = \int_{\partial  D_\alpha } \tilde{\phi}(\tilde{x}) \Scal_{D} \Big(\frac{\partial u}{\partial \nu}\Big)(\tilde{x})
 d\sigma_\alpha (\tilde{x}).
 \end{align*}
Define $\phi := \tilde{\phi} \circ \Psi_\alpha $. Then we get from
\eqref{sigexp} that
 \begin{align*}
 & \int_{\partial  D} \frac{\partial  (\Scal_{D_\alpha } \tilde{\phi})}{\partial  \nu}(x) f(x)
 d\sigma(x) \\
 & = \int_{\partial  D} \phi(x) \Big[ \Scal_{D}
 \Big(\frac{\partial u}{\partial\nu }\Big)(x)+ \alpha  h(x) \frac{\partial }{\partial \nu} \Scal_{D}
 \Big(\frac{\partial u}{\partial \nu}\Big)\Big|_{+}(x) + O(\alpha ^{1+\eta}) \Big] \\
 & \qquad\qquad . (1 -\alpha  \tau(x)h(x)+ O(\alpha ^2)) d\sigma(x) \\
 & = \int_{\partial  D} (\Scal_D \phi) \frac{\partial u}{\partial \nu} d\sigma + \alpha  \int_{\partial D}
 \Big[ -\Scal_D(\tau h \phi) + \big(\frac{1}{2}I + \Kcal_D \big)(h\phi) \Big] \frac{\partial u}{\partial \nu} d\sigma
 + O(\alpha ^{1+\eta}) \\
 & = \int_{\partial  D} \frac{\partial (\Scal_D \phi)}{\partial \nu} \Big|_{-}  f d\sigma + \alpha  \int_{\partial  D}
 \Big[ -\frac{\partial (\Scal_D(\tau h \phi))}{\partial \nu}\Big|_{-}  + \frac{\partial (\Dcal_D(h\phi))}{\partial \nu}  \Big] f d\sigma
 + O(\alpha ^{1+\eta}).
 \end{align*}
Thus we get the following expansion:
 \begin{equation} \label{eq:sexp}
 \frac{\partial  (\Scal_{D_\alpha } \tilde{\phi})}{\partial  \nu}(x) =
 \big(-\frac{1}{2} I + \Kcal_D^* \big) \phi(x) + \alpha  \Lcal_D \phi (x) + O(\alpha ^{1+\eta}),
 \end{equation}
where
 \begin{equation} \label{lform}
 \Lcal_D \phi = \big(\frac{1}{2} I - \Kcal_D^* \big) (\tau h \phi)
 + \frac{\partial (\Dcal_D(h\phi))}{\partial \nu} \quad \mbox{on } \partial D.
 \end{equation}

\subsection{Heat equation}
Recall that $U$ a solution of
\begin{eqnarray}
  (\partial_s -\Delta)U(x,s) &=& 0 \mbox{ in } \Omega\backslash\overline{B}\times (0,S),\label{e1}\\
  (\partial_s -k\Delta)U(x,s) &=& 0 \mbox{ in } B\times (0,S),\label{e2}\\
   \frac{\partial U}{\partial\nu}|_+(x,s) -k \frac{\partial U}{\partial\nu}|_-(x,s) &=& 0 \mbox{ on } \partial B\times (0,S),\label{e3}\\
  U_+(x,s) -U_-(x,s) &=& 0 \mbox{ on } \partial B\times (0,S),\label{e4}\\
   U &=& f \mbox{ on } \partial \Omega\times (0,S).\label{e5}
\end{eqnarray}
Expanding
\begin{equation}\label{eq1}
U(x, s) = U_0(x, s) + \alpha U_1(x, s) + \alpha^2U_2(x, s) + \cdots,\quad x\in \Omega.
\end{equation}
If we substitute the expansions in (\ref{eq1}) into (\ref{e1}), (\ref{e2}), and (\ref{e4}) we formally get for all positive integer $n$, $U_n|_{\partial \Omega\times (0,S)}=f\delta_{0n}$,
$$
(\partial_s -\Delta)U_n(x,s) = 0 \mbox{ in } \Omega\backslash\overline{B}\times (0,S)
$$
and
$$ (\partial_s -k\Delta)U_n(x,s) = 0 \mbox{ in } B\times (0,S).
$$
Let $\widetilde{x} =x + \alpha h(x)\nu(x)$, we have
\begin{equation}\label{eq2}
U(\widetilde{x}, s) = U_0(x, s) +\alpha h(x)\frac{\partial U_0}{\partial\nu}(x,s)+ \alpha U_1(x, s) + \alpha^2U_2(x, s) + \cdots
\end{equation}
If we substitute the expansions in (\ref{eq2}) into (\ref{e4}), we formally get
$$U_0(x, s)|_+ - U_0(x, s)|_-=0;\quad (x,s)\in  \partial B\times (0,S),
$$
and
$$U_1(x, s)|_+ - U_1(x, s)|_-=h(x)(1-k)\frac{\partial U_0}{\partial\nu}(x,s);\quad (x,s)\in  \partial B\times (0,S).
$$
For $\widetilde{x} =x + \alpha h(x)\nu(x)$, we have
\begin{equation}\label{eq3}
\nabla U(\widetilde{x}, s)\cdot\widetilde{\nu}(\widetilde{x}) = \Big(\nabla U_0(x, s) +\alpha h(x) \nabla^2 U_0(x,s)\nu(x)+ \alpha \nabla U_1(x, s) + \cdots \Big)\Big( \nu^{(0)} +\alpha \nu^{(1)}+\cdots \Big).
\end{equation}
If we substitute the expansions in (\ref{eq3}) into (\ref{e2}), we formally get
$$\frac{\partial U_0}{\partial\nu}(x,s)|_+ - k\frac{\partial U_0}{\partial\nu}(x,s)|_-=0;\quad (x,s)\in  \partial B\times (0,S),
$$
and
\begin{equation}\label{eq4}
\frac{\partial U_1}{\partial\nu}(x,s)|_+ - k\frac{\partial U_1}{\partial\nu}(x,s)|_-=(1-k)\frac{d}{dt}\big(h\frac{dU_0}{dt}(x,s)\big);\quad (x,s)\in  \partial B\times (0,S).
\end{equation}
To prove (\ref{eq4}), we expressed the Laplacian in the local coordinates
$$
\Delta=\frac{\partial^2}{\partial \nu^2}-\tau \frac{\partial}{\partial\nu}+(\frac{d}{dt})^2,
$$
then, for $(x,s)\in  \partial B\times (0,S)$
$$
\frac{\partial^2U_0(x,s)}{\partial \nu^2}|_+=\tau\frac{\partial U_0(x,s)}{\partial \nu}|_+-(\frac{d}{dt})^2 U_0(x,s)+\partial_s  U_0(x,s)
$$
and
$$
\frac{\partial^2U_0(x,s)}{\partial \nu^2}|_-=\tau\frac{\partial U_0(x,s)}{\partial \nu}|_- -(\frac{d}{dt})^2 U_0(x,s)+\frac{1}{k}\partial_s  U_0(x,s).
$$
The following theorem holds.
\begin{theorem}\label{thm-heat}
  We formally obtain the following asymptotic formula
\begin{equation}\label{eq5}
U_\alpha(x, s) = U_0(x, s) + \alpha U_1(x, s) + O(\alpha^2); \quad (x,s)\in \Omega\times (0,S),
\end{equation}
where $U_0$ is the unique solution of
$$\left\{
  \begin{array}{ll}
 \displaystyle  (\partial_s -\Delta)U_0(x,s) = 0& \mbox{ in } \Omega\backslash\overline{B}\times (0,S), \\
  \displaystyle (\partial_s -k\Delta)U_0(x,s) = 0 &\mbox{ in } B\times (0,S), \\
    \displaystyle \frac{\partial U_0}{\partial\nu}(x,s)|_+ -k \frac{\partial U_0}{\partial\nu}(x,s)|_- = 0 &\mbox{ on } \partial B\times (0,S),\\
   \displaystyle  U_0(x,s)|_+ -U_0(x,s)|_- = 0 &\mbox{ on } \partial B\times (0,S),\\
\displaystyle U_0(x,0) = \Phi(x)& \mbox{ in } \Omega\times (0,S),\\
  \displaystyle   U_0 = f& \mbox{ on } \partial \Omega\times (0,S).
  \end{array}
\right.
$$
and $U_1$ is the unique solution of
$$\left\{
  \begin{array}{ll}
 \displaystyle  (\partial_s -\Delta)U_1(x,s) = 0& \mbox{ in } \Omega\backslash\overline{B}\times (0,S), \\
  \displaystyle (\partial_s -k\Delta)U_1(x,s) = 0 &\mbox{ in } B\times (0,S), \\
    \displaystyle \frac{\partial U_1}{\partial\nu}(x,s)|_+ -k \frac{\partial U_1}{\partial\nu}(x,s)|_- = (1-k)\frac{d}{dt}\big(h\frac{dU_0}{dt}(x,s)\big) &\mbox{ on } \partial B\times (0,S),\\
   \displaystyle  U_1(x,s)|_+ -U_1(x,s)|_- = h(x)(1-k)\frac{\partial U_0}{\partial\nu}(x,s) &\mbox{ on } \partial B\times (0,S),\\
\displaystyle U_1(x,0) = 0& \mbox{ in } \Omega\times (0,S),\\
  \displaystyle   U_1 = 0& \mbox{ on } \partial \Omega\times (0,S).
  \end{array}
\right.
$$
\end{theorem}

\subsection{Wave equation}
Expanding
\begin{equation}\label{eq6}
u_\alpha(x, s) = u_0(x, s) + \alpha u_1(x, s) + \alpha^2u_2(x, s) + \cdots,\quad x\in \Omega.
\end{equation}
By the same method as above, for the heat equation, we get the following theorem.
\begin{theorem}\label{thm-wave}
  We formally obtain the following asymptotic formula
\begin{equation}\label{eq5}
u(x, s) = u_0(x, s) + \alpha u_1(x, s) + O(\alpha^2); \quad (x,s)\in \Omega\times (0,S),
\end{equation}
where $u_0$ is the unique solution of
$$\left\{
  \begin{array}{ll}
 \displaystyle  (\partial_s^2 -\Delta)u_0(x,s) = 0& \mbox{ in } \Omega\backslash\overline{B}\times (0,S), \\
  \displaystyle (\partial_s^2 -k\Delta)u_0(x,s) = 0 &\mbox{ in } B\times (0,S), \\
    \displaystyle \frac{\partial u_0}{\partial\nu}(x,s)|_+ -k \frac{\partial u_0}{\partial\nu}(x,s)|_- = 0 &\mbox{ on } \partial B\times (0,S),\\
   \displaystyle  u_0(x,s)|_+ -u_0(x,s)|_- = 0 &\mbox{ on } \partial B\times (0,S),\\
\displaystyle u_0(x,0) = \varphi(x),\quad \partial_su_0(x,0) = \psi(x)& \mbox{ in } \Omega\times (0,S),\\
  \displaystyle   u_0 = f& \mbox{ on } \partial \Omega\times (0,S).
  \end{array}
\right.
$$
and $u_1$ is the unique solution of
$$\left\{
  \begin{array}{ll}
 \displaystyle  (\partial_s^2 -\Delta)u_1(x,s) = 0& \mbox{ in } \Omega\backslash\overline{B}\times (0,S), \\
  \displaystyle (\partial_s^2 -k\Delta)u_1(x,s) = 0 &\mbox{ in } B\times (0,S), \\
    \displaystyle \frac{\partial u_1}{\partial\nu}(x,s)|_+ -k \frac{\partial u_1}{\partial\nu}(x,s)|_- = (1-k)\frac{d}{dt}\big(h\frac{du_0}{dt}(x,s)\big) &\mbox{ on } \partial B\times (0,S),\\
   \displaystyle  u_1(x,s)|_+ -u_1(x,s)|_- = h(x)(1-k)\frac{\partial u_0}{\partial\nu}(x,s) &\mbox{ on } \partial B\times (0,S),\\
\displaystyle u_1(x,0) = 0,\quad \partial_su_0(x,0) =0& \mbox{ in } \Omega\times (0,S),\\
  \displaystyle   u_1 = 0& \mbox{ on } \partial \Omega\times (0,S).
  \end{array}
\right.
$$
\end{theorem}

\section{Proofs}
\subsection{Proof of Theorem 4.1}
We formally wrote the asymptotic formula for heat equation in Theorem 4.1 and, in this subsection, we give the proof by using Lemma 3.2. The double-layer operator $ \Dcal_D$ has the integral representation
$$ (\Dcal_D \psi) (x,t) =-\frac{1}{2\pi}
 \int_{0}^{t}\int_{\partial  D}\frac{\psi(x,\tau)}{4(t-\tau)^2} \nu_y\cdot (y-x)\exp(-\frac{|x-y|^2}{4(t-\tau)})d \sigma (y)d\tau.
$$

To prove the main result, we return
to the system of integral equations \eqref{eq:syst1} and derive
asymptotic expansions of $u_\alpha  $ as $\alpha  \to 0$, considering
separately each of the three cases.

Recall from \cite{AK2} that
 $$
 ||\partial^i H_\alpha  ||_{\mathcal{C}^l(\Omega^\prime)} \leq C, \quad i \in \mathbb{N}^2, l \in
 \mathbb{N},
 $$
for any $\Omega^\prime \subset \subset \Omega$, where $C$ depends
only on $\mbox{dist}(\Omega^\prime,
\partial \Omega), \Omega, g, i,$ and $l$. Then the Taylor expansion of $\partial H_\alpha /
\partial \tilde{\nu}  (\tilde{x})$ writes:
\begin{equation}\label{eq:ha}
\frac{\partial H_\alpha }{\partial \tilde \nu} (\tilde{x})=\frac{\partial H_\alpha }{\partial
\nu}(x)+\alpha  \Big (-h'(t)\frac{\partial H_\alpha }{\partial T}(x)+h(t)\langle D^2
H_\alpha (x)\nu(x),\nu(x)\rangle \Big )+ O(\alpha ^2),
\end{equation}
for $\tilde{x} = \Psi_\alpha (x)$, $x \in \partial D$.

Using Lemma \ref{lemkd} together with equations \eqref{eq:stilde},
\eqref{eq:sexp}, and \eqref{eq:ha}, we can write \eqref{eq:syst1}
under the following form: for $x \in \partial D$,
\begin{equation} \label{eq:un} \left\{
\begin{array}{l}
\displaystyle (\lambda I-\Kcal_{D}^*) {\phi_\alpha }(x)- (\frac{1}{2}
I+\Kcal_{D}^*)\psi_\alpha (x)- \alpha  \big[ \Kcal_D^{(1)} \phi_\alpha  +
\Rcal_D\psi_\alpha (x) \big]+ O(\alpha ^{1+\eta}) \\
\nm \quad \quad \displaystyle =\frac{\partial H_\alpha }{\partial \nu}(x) + \alpha  \Big
(-h'(t)\frac{\partial H_\alpha }{\partial T}(x)+h(t)\langle
D^2 H_\alpha (x)\nu(x),\nu(x)\rangle \Big )+O(\alpha^2), \\
\nm \displaystyle (\lambda I+\Kcal_{D}^*)\psi_\alpha (x)+(-\frac{1}{2}
I+\Kcal_{D}^*) {\phi_\alpha }(x)+\alpha \Lcal_D
{\phi_{\alpha }}(x)+O(\alpha^{1+\eta})=-\frac{\partial H_\alpha }{\partial \nu}(x),
\end{array} \right.
\end{equation}
where $\phi_\alpha  = \tilde{\phi_\alpha }\circ \Psi_\alpha $.
\subsection{Proof of Theorem 4.2}
As we said above, we have formally wrote the asymptotic formula for heat equation in Theorem 4.2 and, in this subsection, we give the proof by using Lemma 3.2. 

\section{The identification procedure}
Our dynamic identification procedure can be regarded as a first attempt to generalize that
of [1, 2, 6] to the time-dependent case. It yields simple and effective computational identification algorithms.\\

In this section, we give a reconstruction formula of the shape deformation h from dynamical
measurements of boundary perturbations $U_\alpha-U_0$ for heat conduction problem (resp. $u_\alpha-u_0$ for wave equation) on $\partial\Omega$.

\subsection{Reconstruction of the interface deformation for heat equation}
Recall from Theorem \ref{thm-heat} that
\begin{equation}\label{eq7}
U_\alpha(x, s) = U_0(x, s) + \alpha U_1(x, s) + O(\alpha^2); \quad (x,s)\in \Omega\times (0,S).
\end{equation}
Then, we can assume that we are in possession of the boundary measurements of
$$U_\alpha(x, s)\quad \mbox{on }\partial\Omega\times (0,S)
$$
for a particular choice of data $\phi$ and $f$ which implies that the background solution $U$ of
the heat equation (\ref{main-heat0}) in the absence of deformation of background medium can be given
explicitly in terms of Fourier-Bessel series. Let $V$ be the solution of the following problem
\begin{equation}\label{main-heat0-adj}
\left\{
  \begin{array}{ll}
    \partial_s V-\nabla\cdot \big(k^\star\nabla V\big) = 0 & \mbox{ in } \Omega\times (0,S), \\
    V(x,s=0)=\psi(x)&\mbox{ in } \Omega,\\
   \frac{\partial V}{\partial\nu}=g &\mbox{ on }\partial\Omega\times(0,S),\\
    V=g_1&\mbox{ on }\partial\Omega\times(0,S),
  \end{array}
\right.
\end{equation}
Here $k^\star=\chi_{\Omega\setminus \overline{B}}+k\chi_{B}$. If we multiply $(U_\alpha-U_0)$ by the data $g$ and if we integrate by parts
over $\partial\Omega\times(0, S)$, we immediately get:
\begin{equation}\label{eq8}
  \int_{0}^{S}\int_{\partial\Omega}g(U_\alpha-U_0)d\sigma ds=\alpha \int_{0}^{S}\int_{\partial\Omega}gU_1d\sigma ds+O(\alpha^2)=\alpha \int_{0}^{S}\int_{\partial\Omega}g\frac{\partial V}{\partial\nu}d\sigma ds+O(\alpha^2)
\end{equation}
$$=\alpha \frac{1-k}{k}\int_{0}^{S}\int_{\partial B}h\Big[k\frac{\partial U}{\partial\nu}|_-\frac{\partial V}{\partial\nu}|_- -\frac{\partial V}{\partial T}\frac{\partial U}{\partial T} \Big]d\sigma ds+O(\alpha^2),
$$
where $U_0 = U$ and $U_\alpha$ are solutions of problems (\ref{main-heat0}) and (\ref{u-heat-pert}) respectively.
Now, our goal is to use the formula (\ref{eq8}) to determine an approximation of $h$ by appropriately
choosing $V$. This procedure of construction is based on the extraction of Fourier coefficients
$h_p$, related to function $h$, from a finite number of measurements provided that the order of
magnitude of $h_p$ is much larger than $\alpha$.\\

To proceed with this procedure of reconstruction, we fix $\Omega$ to be the unit disk centered at
the origin, and $B$ to be a two-dimensional disk centered at the origin with radius $a$. Set the
initial temperature distribution $\Phi(r,\theta)= r^{|m|}e^{im\theta}$, then it follows that
$$U(r,\theta; s) = A_{mn}J_{|m|}(\mu_{mn}r/a)e^{im\theta}e^{-k\mu_{mn}^2s/a^2}\quad \mbox{for } r \leq a,
$$
where,
$$A_{mn}:=\frac{2a^m}{\mu_{mn}J_{m+1}(\mu_{mn})},\quad m\geq 0,
$$
and where $\mu_{mn}$ is the n-th positive zero of the Bessel function of the first kind $J_m(x)$ (see
for example \cite[pp. 37-39]{bowman}.\\
Set, $g(\theta, s) =e^{-ks(\mu_{pq}/a)^2}e^{ip\theta}$, and  $\psi(r,\theta) = r^{|p|}e^{ip\theta}$. Then we can get:
$$V(r,\theta; s) = B_{pq}J_{|p|}(\mu_{pq}r/a)e^{ip\theta}e^{-ks(\mu_{pq}/a)^2}\quad \mbox{for } r \leq a,
$$
where,
$$B_{pq}:=\frac{2a^p}{\mu_{pq}J_{p+1}(\mu_{pq})},\quad p\geq 0.
$$
It follows from (\ref{eq8}) that
 \begin{equation}\label{eq9}
  \int_{0}^{S}\int_{\partial\Omega}g(U_\alpha-U_0)d\sigma ds=\alpha C_{mp}\int_{0}^{2\pi}he^{i(m+p)\theta}d\theta+O(\alpha^2);
 \end{equation}
where the constant
$$C_{mp}:=a\frac{1-k}{k}\sum_{n,q\geq1}A_{mn}B_{pq}\Big( \mu_{mn}\mu_{pq}J_{m}^{\prime}(\mu_{mn})J_{p}^{\prime}(\mu_{pq})+\frac{mp}{k}J_{m}(\mu_{mn})J_{p}(\mu_{pq})\Big)\frac{1
-e^{-Sk(\frac{\mu_{mn}^2+\mu_{pq}^2}{a^2})}}{\mu_{mn}^2+\mu_{pq}^2}.
$$
Relation (\ref{eq9}) implies that the Fourier coefficient $h_j$ can be determined from measurements
of $U_\alpha-U$ on $\partial\Omega\times(0, S)$ provided that the order of magnitude of $|h_j|$ is much larger than $\alpha$. To reconstruct higher Fourier coefficients of $h$ or more accurately the first ones, the
high-order asymptotic expansions derived in this paper should be used.
\subsubsection{Algorithm}
The scheme is relatively simple. We use the leading-order term to determine the first Fourier
coefficients of $h$, say $h_j$, for $|j|\leq N$. Then, we insert these coefficients into the second-order
term by approximating $h$ by $\sum_{|j|\leq N}h_je^{ij\theta}$ and we apply the asymptotic expansion in Theorem
\ref{thm-heat} to both determine further Fourier coefficients and get more accurate approximations
of $(h_j)_{|j|\leq N}$.\\

On the other hand, the arclength parametrization of the boundary of the disk $B$ can
be given by $X(\theta) = a(\cos \theta, \sin \theta)$. Assume that $h(\theta) =\sum_{|j|\leq N}h_je^{ij\theta}$, where $N$ is a given
positive integer. If we restrict Fourier series of $h$ in relation (\ref{eq9}) to order $N$, we immediately
get
 \begin{equation}\label{eq10}
  \int_{0}^{S}\int_{\partial\Omega}g(U_\alpha-U_0)d\sigma ds=\alpha(\frac{1}{2\pi})\sum_{|j|\leq N}\sum_{m\geq 0} C_{m -j-m}\int_{0}^{2\pi}h_j+O(\alpha^2).
 \end{equation}

This method can be repeated recursively. Let us neglect the asymptotically small remainder
in the asymptotic formula (\ref{eq10}), then we have the following scheme:
 \begin{equation}\label{eq11}
\Lambda(\alpha):=  \int_{0}^{S}\int_{\partial\Omega}g(U_\alpha-U_0)d\sigma ds=z_{-N}h_{-N} +\cdots + z_0h_0 + z_1h_1 +\cdots + z_Nh_N;
 \end{equation}
where
$$z_j =\frac{\alpha}{2\pi}\sum_{m\geq 0} C_{m -j-m},\quad    \mbox{for } j =-N,\cdots, N.
$$

To determine the $(2N + 1)$ coefficients $h_{-N},\cdots,h_0,\cdots, h_N$ we need to define a number of
$2N$-other linear equations where these coefficients are their unknown. Let $j_0$ be the smallest
positive integer such that $h_p\ll \frac{1}{10^{j_0}}$, for $|p|\leq N$. Then, by replacing $\alpha$ by $\frac{1}{10^{j_0}}$ and by
taking into account relation (\ref{eq11}), we obtain
$$
\Lambda_{j_0}=\Lambda(\frac{1}{10^{j_0}})=z_{-N}(j_0)h_{-N} +\cdots + z_0(j_0)h_0 + z_1(j_0)h_1 +\cdots + z_N(j_0)h_N.
$$
Thus, taking values greater than $j_0$, one can find a sufficient number of linear equations to
find the Fourier coefficients $h_j$ and therefore the function $h$ is approximately reconstructed.
\subsection{Reconstruction of the interface deformation for wave equation}
Firstly, we recall from Theorem \ref{thm-wave} that
\begin{equation}\label{eq12}
u_\alpha(x, s) - u_0(x, s)= \alpha u_1(x, s) + O(\alpha^2); \quad (x,s)\in \Omega\times (0,S).
\end{equation}
As done for reconstruction via dynamical boundary measurements for heat equation, we
firstly introduce the function $v$ as solution of the following initial boundary value problem
for the (scalar) wave equation
\begin{equation}\label{main-heat0-adj}
\left\{
  \begin{array}{ll}
    \big(\partial_s^2 -\nabla\cdot c^\star\nabla \big)v = 0 & \mbox{ in } \Omega\times (0,S), \\
    v(x,s=0)=v^0(x),\quad  \partial_s v(x,s=0)=v^1(x) &\mbox{ in } \Omega,\\
   \frac{\partial v}{\partial\nu}=g &\mbox{ on }\partial\Omega\times(0,S),\\
    v=g_1&\mbox{ on }\partial\Omega\times(0,S),
  \end{array}
\right.
\end{equation}
Here $c^\star=\chi_{\Omega\setminus \overline{B}}+c\chi_{B}$. Multiplying (\ref{eq12}) by the data $g$ and integrating by parts over $\partial\Omega\times(0,S)$, we immediately get:
\begin{equation}\label{eq13}
  \int_{0}^{S}\int_{\partial\Omega}g(u_\alpha-u_0)d\sigma ds=\alpha \int_{0}^{S}\int_{\partial\Omega}gu_1d\sigma ds+O(\alpha^2)=\alpha \int_{0}^{S}\int_{\partial\Omega}g\frac{\partial v}{\partial\nu}d\sigma ds+O(\alpha^2)
\end{equation}
$$=\alpha \frac{1-c}{c}\int_{0}^{S}\int_{\partial B}h\Big[ -\frac{\partial v}{\partial T}\frac{\partial u}{\partial T}+c\frac{\partial u}{\partial\nu}|_-\frac{\partial v}{\partial\nu}|_- \Big]d\sigma ds+O(\alpha^2),
$$
where $u_0 = u$ and $u_\alpha$ are solutions of problems (\ref{main-wave0}) and (\ref{u-wave-pert}) respectively.\\

The objective, now, is to formulate an approximation of $h$ by using relation (\ref{eq12}) and by
appropriately choosing the function $v$. As stated in the previous section, this procedure of
identification is based on the extraction of Fourier coefficients $h_n$, associated to function $h$,
from a finite number of measurements provided that the order of magnitude of $h_n$ is much
larger than $\alpha$.\\
To illustrate this procedure of reconstruction, we suppose that $\Omega$ is the unit disk centered
at the origin, and we assume that $B$ is a two-dimensional disk centered at the origin with
radius $b$. Choosing the initial wave distributions $u^0(r,\theta)= r^{|m|}e^{im\theta}$, $u^1(r,\theta)= \varrho(r) e^{im\theta}$ (where $\varrho(r)$ is a radial function). Then, it follows that
$$u(r,\theta; s) = D_{mn}J_{|m|}(\lambda_{mn}r/b)e^{im\theta}e^{-c\lambda_{mn}^2s/b^2}\quad \mbox{for } r \leq b,\quad m\geq 0,\quad n\geq 1.
$$
where $\lambda_{mn}$ is the $n$-th positive zero of the Bessel function $J_m(x)$, and $D_{mn}$ are constants
which could be calculated as integrals by using the orthogonality relationship for Bessel
functions combined with initial conditions.\\

Set, $g(\theta, s) =e^{ics\lambda_{pq}/b}e^{ip\theta}$, and  $v^0(r,\theta) = r^{|p|}e^{ip\theta}$, $v^1(r,\theta) = \varrho(r)e^{ip\theta}$. Then we can get:
$$v(r,\theta; s) = E_{pq}J_{|p|}(\lambda_{pq}r/b)e^{ip\theta}e^{ics\lambda_{pq}/b}\quad \mbox{for } r \leq b,
$$
where the coefficients $E_{pq}$ could be calculated as done for $D_{mn}$. It follows from (\ref{eq12}) that
 \begin{equation}\label{eq14}
  \int_{0}^{S}\int_{\partial\Omega}g(u_\alpha-u_0)d\sigma ds=\alpha F_{mp}\int_{0}^{2\pi}he^{i(m+p)\theta}d\theta+O(\alpha^2);
 \end{equation}
where the constant
$$F_{mp}:=\frac{1-c}{c}\sum_{n,q\geq1}iD_{mn}E_{pq}\Big( \lambda_{mn}\lambda_{pq}J_{m}^{\prime}(\lambda_{mn})J_{p}^{\prime}(\lambda_{pq})+\frac{mp}{c}J_{m}(\lambda_{mn})J_{p}(\lambda_{pq})\Big)\frac{1
-e^{icS(\frac{\lambda_{mn}+\lambda_{pq}}{b})}}{\lambda_{mn}+\lambda_{pq}}.
$$
Relation (\ref{eq14}) implies that the Fourier coefficient $h_j$ can be determined from measurements
of $u_\alpha-u$ on $\partial\Omega\times(0, S)$ provided that the order of magnitude of $|h_j|$ is much larger than $\alpha$. To reconstruct higher Fourier coefficients of $h$ or more accurately the first ones, the high-order asymptotic expansions derived in this paper should be used.\\


\section*{Acknowledgments}
This work is partially supported by Institute for Advanced Studies and AGM Research Center in Mathematics (UMR CNRS 8088), CY Cergy Paris Universit\'e, F-95302 France.


\end{document}